\newcommand{\mathbb}[1]{{\bf #1}}

\documentclass[fancyheadings,twoside,psfig,float,fleqn,epsf,12pt]{article}
\usepackage [dvips]{graphicx}
\usepackage{psfig}
\usepackage{epsf}
\usepackage{amssymb}
\usepackage{float}
\usepackage{fancyheadings}

\restylefloat{figure}

\newcommand{\be}{\begin{equation}}
\newcommand{\ee}{\end{equation}}
\newcommand{\bea}{\begin{eqnarray}}
\newcommand{\eea}{\end{eqnarray}}
\newcommand{\bean}{\begin{eqnarray*}}
\newcommand{\eean}{\end{eqnarray*}}

\newcommand{\pp}{\prime}
\newcommand{\LL}{\mbox{\tiny L}}
\newcommand{\RR}{\mbox{\tiny R}}
\newcommand{\CC}{\mbox{\tiny C}}

\newcommand{\uj}{\bar{u}_{j}^{n}}
\newcommand{\ujp}{\bar{u}_{j+1}^{n}}
\newcommand{\ujm}{\bar{u}_{j-1}^{n}}

\newenvironment{remark}{{\flushleft \bf Remark:}}{}
\newenvironment{remarks}{{\flushleft \bf Remarks:}}{}

\newenvironment{acknowledgment}{{\flushleft \bf Acknowledgment:}}{}

\begin{document}

\title{A Third-Order Semi-Discrete Central Scheme\\ 
       for Conservation Laws and Convection-Diffusion Equations}
\author{Alexander Kurganov\footnotemark[2] \and Doron Levy\footnotemark[3]}
\date{}

\renewcommand{\thefootnote}{\fnsymbol{footnote}}
\footnotetext[2]
 {Department of Mathematics, University of Michigan, Ann Arbor,
  MI 48109; {\tt kurganov@math.lsa.umich.edu}}
\footnotetext[3]
 {Department of Mathematics, University of California, Berkeley, CA 94720, 
  and Lawrence Berkeley National Lab; {\tt dlevy@math.berkeley.edu}}
\renewcommand{\thefootnote}{\arabic{footnote}}

\maketitle


\begin{abstract}
We present a new third-order, semi-discrete, central method  for 
approximating solutions to multi-dimensional
systems of hyperbolic conservation laws, convection-diffusion
equations, and related problems.  Our method is a high-order 
extension of the recently proposed second-order,
semi-discrete method in \cite{kurganov-tadmor:semi}.

The method is derived independently of the specific piecewise polynomial
reconstruction which is based on the previously computed cell-averages.  
We demonstrate our results, by focusing on the new third-order CWENO
reconstruction presented in \cite{levy-puppo-russo:balance}.
The numerical results we present, show the desired accuracy, 
high resolution and robustness of our method.
\end{abstract}


\bigskip
\noindent
{\bf AMS(MOS) subject classification.} Primary 65M10; secondary 65M05.



\section{Introduction}                     \label{section:introduction}
\setcounter{equation}{0}
\setcounter{figure}{0}
\setcounter{table}{0}

Numerical methods for approximating solutions of hyperbolic
conservation laws,
\be
\frac{\partial}{\partial t}u(x,t)+\frac{\partial}{\partial x}f(u(x,t))=0,
\label{CL}
\ee
and of the related convection-diffusion equations,
\be
\frac{\partial}{\partial t}u(x,t)+\frac{\partial}{\partial x}f(u(x,t))=
\frac{\partial}{\partial x}Q[u(x,t),u_x(x,t)],
\label{CD}
\ee
have attracted a lot 
of attention in recent years (see, e.g.,
\cite{god-rav:difference, tadmor:approximate}, 
and the references therein).
Here, $u(x,t)$ is a conserved quantity, $f(u)$ is a nonlinear convection
flux and $Q(u,u_x)$ is a dissipation flux satisfying the (weak)
parabolicity condition,
$\frac{\partial}{\partial s}Q(u,s)\geq 0$, $\forall u,s$.
In the most general case
$u=(u_1,...,u_n)$ is an $n$-vector in the $d$-spatial variables,
$x=(x_1,...,x_d)$, and $f$ and $Q$ are vector-functions.

In this paper, we focus on the class of {\em central schemes\/},
all of which can be viewed as an extension of the well-known
Lax-Friedrichs (LxF) scheme, \cite{lxf}. The first-order LxF method
enjoys the major advantage of simplicity over the 
upwind schemes (e.g., the Godunov scheme, \cite{God}): no (approximate) Riemann 
solvers or characteristic decompositions are involved in its construction,
and therefore, its realization for complicated multi-dimensional
systems is rather simple.
At the same time, the LxF scheme suffers from excessive
numerical dissipation, which causes a poor (smeared) resolution of
discontinuities and rarefaction waves.

A second-order, non-oscillatory central scheme was first introduces
by Nessyahu and Tadmor in \cite{nessyahu-tadmor:non-oscillatory}.
Since then, the Nessyahu-Tadmor (NT)
scheme was further extended to higher orders of
accuracy, \cite{liu-tadmor:3rd} (also see \cite{bianco-puppo-russo:central,
levy-puppo-russo:1d}), as well as to the multi-dimensional systems 
(\ref{CL}), in \cite{arminjon:nt.french} and \cite{jiang-tadmor:nonosc},
(also \cite{levy:third,levy-puppo-russo:2d-3,
levy-puppo-russo:balance, levy-puppo-russo:2d-4}).

The main ingredient in the construction of the NT method is a second-order,
non-oscillatory, MUSCL-type \cite{vLeV}, piecewise linear
interpolant (instead of the piecewise constant one, employed in the LxF scheme)
in combination with the exact solver for the time evolution.
This approach allows to significantly improve the resolution of non-smooth
solutions to hyperbolic conservation laws, (\ref{CL}), while retaining
the main advantage of the LxF scheme -- {\em simplicity}.

Unfortunately, applying the fully-discrete NT scheme (or its higher-order
extensions) to the second-order convection-diffusion equations,
(\ref{CD}), does not provide the desired resolution of discontinuities
(see, e.g., \cite{KLR, KR, kurganov-tadmor:semi}). This loss of resolution
occurs due to the accumulation of excessive numerical dissipation, which
is typical of fully-discrete central schemes with small time-steps
$\Delta t \sim {(\Delta x)}^2$ (see \cite{kurganov-tadmor:semi} for
details).

To circumvent this difficulty, a second-order
{\em semi-discrete\/} central scheme was introduced
by Kurganov and Tadmor in \cite{kurganov-tadmor:semi}.
This scheme has smaller dissipation than the NT scheme, and
unlike the fully-discrete central schemes, it can be efficiently used with
time-steps as small as required by the CFL stability restriction.

The basic idea in the construction of the second-order semi-discrete scheme 
was to use a more accurate information about the 
{\em local\/} speed of propagation of the discontinuities.   
One was then able to derive a non-staggered semi-discrete central method, by
first integrating over non-equally spaced control volumes, out of which
a new piecewise linear interpolant was reconstructed and finally projected
on its cell-averages (without evolving in time).  The final step, was
first introduced in \cite{jiang-levy-osher-tadmor:stg}, in a somewhat different
context of transforming staggered methods into non-staggered methods.

In this paper we extend the results of \cite{kurganov-tadmor:semi} by
introducing a new {\em third-order, semi-discrete, central scheme}.
Our new scheme is derived in a general form which is independent of 
the reconstruction step, as long as the reconstructed interpolant
is sufficiently accurate and non-oscillatory. 
In particular, we use the new third-order CWENO reconstruction proposed
in \cite{levy-puppo-russo:balance}.  This reconstruction provides
a third-order accurate interpolant which is built from the given cell-averages
such that it is non-oscillatory in the essentially non-oscillatory
(ENO) sense (see \cite{harten-eoc:eno}, \cite{shu:eno}). This 
interpolant is written as a convex combination of two one-sided linear
functions and one centered parabola.  In smooth regions this convex
combination guarantees the desired third-order accuracy. It automatically
switches to a second-order, one-sided, linear reconstruction in the 
presence of large gradients.  Such weighted essentially non-oscillatory
(WENO) reconstructions were first introduced in the upwind framework,
\cite{jiang-shu:weno, liu-osher-chan:weno}, after which
they were extended to the central framework, \cite{levy-puppo-russo:1d,
levy-puppo-russo:2d-3, levy-puppo-russo:balance, levy-puppo-russo:2d-4}.

This paper is organized as follows.  We start in \S\ref{section:central}
with a brief overview of central schemes for conservation laws.  In 
particular we focus in \S\ref{subsection:cweno} on the CWENO reconstruction
which we use as the building block for our third-order method below.

We then proceed to construct our new third-order scheme.  
First, we deal with the fully-discrete, one-dimensional setup in
\S\ref{section:fully-discrete}.  
This new fully-discrete scheme is sketched in equation (\ref{eq:fully.scheme}).  
We only give the required details that are necessary to fulfill our final
goal, namely, to derive the semi-discrete scheme.

With the fully-discrete scheme, (\ref{eq:fully.scheme}), we are ready to 
approach the semi-discrete limit in \S\ref{subsec:1D}.
Our new third-order, one-dimensional,
semi-discrete scheme is then summarized in 
equation (\ref{semidiscrete_scheme}).
This scheme is written in a general form which is 
independent of the reconstruction step and can also be combined with
any appropriate ODE solver for carrying out the time evolution. 
In \S\ref{subsec:2D} we then extend our semi-discrete scheme to 
multidimensional hyperbolic and (degenerate) parabolic problems.

We end by presenting several numerical examples in \S\ref{section:numerical},
in which we approximate solutions to hyperbolic conservation laws as well
as to convection-diffusion equations.  Our new method is shown to enjoy
the expected high-accuracy as well as the robustness and the simplicity
of the entire family of central schemes.


\section{Central Schemes for Conservation Laws}   \label{section:central}
\setcounter{equation}{0}
\setcounter{figure}{0}
\setcounter{table}{0}

We briefly overview the framework of central schemes for conservation laws.
Consider the one-dimensional system (\ref{CL}).
To approximate its solutions, we introduce a spatial 
scale, $\Delta x$, and integrate over the cell
 $I(x) := \{ \xi \mid~ |\xi-x| \leq \Delta x / 2 \}$,
\be
\bar{u}_{t}+\frac{1}{\Delta x}\left[f\Bigl(u(x+\frac{\Delta x}{2},t)\Bigr)
           +f\Bigl(u(x-\frac{\Delta x}{2},t)\Bigr)\right] = 0.
 \label{eq:exact.integration}
\ee
Here and below, $\bar{u}$ denotes the average of $u$ over $I$,
\[
\bar{u}(x,t):=\frac{1}{\Delta x}\int\limits_{I(x)} u(\xi,t)d\xi.
\]

Introducing a time scale, $\Delta t$, integrating
in time from $t$ to $t+\Delta t$ and sampling (\ref{eq:exact.integration})
at the cells $[x_j,x_{j+1}]$, we obtain
\be
\bar{u}^{n+1}_{j+1/2} = \bar{u}^{n}_{j+1/2} -
\frac{1}{\Delta x} \int\limits_{\tau = t^{n}}^{t^{n+1}}
\left[ f(u(x_{j+1},\tau)) - f(u(x_j,\tau))
\right] d \tau,
\label{eq:central}
\ee
where $x_{j}:=j \Delta x$, $t^{n}:=n \Delta t$,
$u_{j}^{n}:=u(x_{j}, t^{n})$ and $\uj:=\bar{u}(x_j,t^n)$.
Assuming that at time $t=t^{n}$ we have computed the cell-averages of the
approximate solution, $\{\uj\}$, we would like
to utilize (\ref{eq:central})
to compute the cell-averages at the next time level, $t^{n+1}=t^n+\Delta t$.
To that extent, we introduce a piecewise-polynomial reconstruction,
\be
u(x,t^{n}) \approx \sum_{j} P_{j}(x) \chi_{j}(x),
\label{eq:reconstruction}
\ee
where $\chi_{j}(x)$ is the characteristic function of the cell
$I_{j}:=I(x_j)$, and $P_{j}(x)$ is a polynomial which is reconstructed 
from the computed cell-averages, $\{\uj\}$.
The degree of the polynomial depends on
the desired order of accuracy of the method.
Having such an approximation to $u(x,t^n)$, (\ref{eq:reconstruction}),
we can easily compute the RHS of (\ref{eq:central}).
The first term, $\bar{u}^{n}_{j+1/2}$, equals
\[
\bar{u}^{n}_{j+1/2}=\int\limits_{x_j}^{x_{j+1/2}} P_{j}(x) dx +
\int\limits_{x_{j+1/2}}^{x_{j+1}} P_{j+1}(x) dx.
\]
For a sufficiently small time-step, $\Delta t$, the solution
of (\ref{CL}) subject to the initial data (\ref{eq:reconstruction}),
prescribed at time $t=t^n$, will remain smooth at some neighborhood of
the grid points $x_{j}$ for $t\in[t^n,t^{n+1}]$.
Hence, the integrals on the RHS of
(\ref{eq:central}) can be approximated
using a sufficiently accurate quadrature,
which is determined by the overall desired accuracy of the method.  
The values at the intermediate times which will be required
in the quadrature, can be predicted either by a Taylor expansion or using a 
Runge-Kutta method (consult \cite{bianco-puppo-russo:central,
levy-puppo-russo:1d,
liu-tadmor:3rd, nessyahu-tadmor:non-oscillatory}).

For example, a piecewise-constant reconstruction,
$P_{j}(x) = \bar{u}_{j}^{n}$, and a first-order
quadrature,
\[ \int\limits_{t^{n}}^{t^{n+1}} f(u(t)) dt \sim \Delta t f(\bar{u}^{n}),
\]
will result with the staggered-LxF scheme
(with $\lambda := \Delta t / \Delta x$ denoting the mesh ratio),
\[
\bar{u}_{j+1/2}^{n+1} = \frac{\ujp+\uj}{2}-
\lambda ( f(\bar{u}_{j+1}^{n}) - f(\bar{u}_{j}^{n})).
\]
A piecewise linear reconstruction,
$P_{j}(x) = \bar{u}_{j}^{n} + (u_{x})^{n}_{j}(x-x_{j})$,
with a second-order quadrature in time (such as the mid-point rule), 
results with the Nessyahu-Tadmor (NT) scheme.  Applying nonlinear limiters
on the discrete slopes, $(u_{x})^{n}_{j}$, will prevent oscillations
(for details, see \cite{nessyahu-tadmor:non-oscillatory}).

To obtain a third-order central scheme, one should use a third-order,
piecewise parabolic reconstruction together with a more accurate
quadrature in time, e.g., Simpson's quadrature rule (see
\cite{liu-tadmor:3rd} for details).

\begin{remarks}
\begin{enumerate}
\item {\bf robustness.} In order to reconstruct a non-oscillatory interpolant,
one typically is required to use nonlinear limiters. These limiters decrease
the order of accuracy of the method at extrema and by that they
play a stabilizing role (e.g., see \cite{vLeV,
liu-tadmor:3rd, nessyahu-tadmor:non-oscillatory, tadmor:approximate}).
\item {\bf numerical dissipation and time step.}
When using fully-discrete central schemes to approximate
solutions of convection-diffusion equations, (\ref{CD}),
the stability restriction enforce small time-steps,
$\Delta t \sim {(\Delta x)}^2$. That is why the numerical dissipation is
accumulated and
we do not obtain high resolution of discontinuities
(see \cite{kurganov-tadmor:semi} for details).

This problem can be avoided by using semi-discrete schemes instead of the
fully-discrete schemes.  Such a second-order, central, semi-discrete scheme was
introduced in \cite{kurganov-tadmor:semi}.
In this paper we develop a third-order, central, semi-discrete scheme with
small numerical dissipation, which
can be efficiently used with the small time-steps
required due to the second-order operators.
\item {\bf upwind schemes.}
Sampling (\ref{eq:exact.integration}) at the cells $I_j$,
will result with upwind schemes.  
Here, one remains with the discontinuities along
the interfaces and is bound to solve the Riemann problems there, or at least to 
approximate their solutions.
In the scalar, one-dimensional case this can be easily accomplished,
but the Riemann problem has no known solution in the general case of 
systems and/or several space dimensions.

This is the reason for why central schemes can be considered as universal methods for
solving hyperbolic conservation laws: Riemann solvers are not
involved in their construction, and moreover, since (\ref{eq:central}) can be
carried out componentwise, no characteristic decomposition
is required.
\end{enumerate}
\end{remarks}


\subsection{CWENO reconstruction}          \label{subsection:cweno}

The first one-dimensional, third-order central scheme in \cite{liu-tadmor:3rd},
implemented the non-oscillatory piecewise parabolic reconstruction
proposed by Liu and Osher in \cite{liu-osher:nonosc}.
Since then, a variety of simpler reconstructions has appeared in the literature.
Among these, we would like to mention the Central-ENO
reconstruction in \cite{bianco-puppo-russo:central} and the
Central-WENO (CWENO) reconstruction in
\cite{levy-puppo-russo:1d} and \cite{levy-puppo-russo:balance},
which was extended
to the two-dimensional setup in \cite{levy-puppo-russo:2d-3} and
\cite{levy-puppo-russo:2d-4}. 

Our new third-order semi-discrete method which we develop in 
\S\ref{section:fully-discrete} and \S\ref{section:semi-discrete} below, can be
integrated with any third-order, non-oscillatory reconstruction.
In our numerical simulations presented in \S\ref{section:numerical}, 
we will use the method recently presented in \cite{levy-puppo-russo:balance}, 
which we will now briefly overview.

In each cell $I_{j}$ we reconstruct
a quadratic polynomial as a convex combination
of three polynomials $P_{\LL}(x), P_{\RR}(x)$ and $P_{\CC}(x)$, 
\be
 P_{j}(x) = w_{\LL} P_{\LL}(x) + w_{\RR} P_{\RR}(x) + w_{\CC} P_{\CC}(x),
 \label{eq:cweno.reconstrcution}
\ee
with positive weights $w_{i} \geq 0, \forall i\in\{\CC,\RR,\LL\}$, and
 $\sum_{i} w_{i} = 1$.  The polynomials $P_{\LL}(x), P_{\RR}(x)$
correspond to left and right one-sided linear reconstructions, respectively, 
while $P_{\CC}(x)$ is a parabola, centered around $x_{j}$.

The linear functions, $P_{\RR}(x)$ and $P_{\LL}(x)$, are uniquely determined
by requiring them to conserve the one-sided cell averages 
($\uj, \ujp$ and $\uj, \ujm$, respectively) as
\be
P_{\RR}(x) = \uj + \frac{\ujp - \uj}{\Delta x}(x-x_{j}),
\quad
P_{\LL}(x) = \uj + \frac{\uj - \ujm}{\Delta x}(x-x_{j}).  \label{eq:cweno.side}
\ee

The centered parabola, $P_{\CC}(x)$, is chosen such as to satisfy,
\be
P_{\mbox{\tiny EXACT}}(x)=
c_{\LL} P_{\LL}(x) + c_{\RR} P_{\RR}(x) + (1-c_{\LL}-c_{\RR}) P_{\CC}(x),
\label{P_exact}
\ee
with constants $c_{i}$'s.
Here, $P_{\mbox{\tiny EXACT}}(x)$ is the unique parabola that conserves
the three cell averages, $\ujm, \uj$ and $\ujp$, which is given by
\be
P_{\mbox{\tiny EXACT}}(x)=
u_{j}^{n} + u_{j}^{\pp}(x-x_{j}) + \frac{1}{2} u_{j}^{\pp\pp}(x-x{j})^{2}.
\label{eq:cweno.exact}
\ee
The approximations to the point-values of $u(x_j,t^n), u_x(x_j,t^n)$ and
$u_{xx}(x_j,t^n)$, are denoted by $u_{j}^{n}, u_{j}^{\pp}, u_{j}^{\pp \pp}$ and
are given by
\bean
u_{j}^{n} & = & \uj - \frac{1}{24}(\ujp-2\uj + \ujm), \\
\qquad u_{j}^{\pp} & = & \frac{\ujp - \ujm}{2 \Delta x}, \qquad
u_{j}^{\pp \pp}   =  \frac{\ujm - 2\uj + \ujp}{\Delta x^2}
\eean
In \cite{levy-puppo-russo:balance} it was shown that every symmetric selection
of the constants $c_{i}$'s in (\ref{P_exact}) will provide the desired 
third-order accuracy.
For example, by taking, $c_{\LL} \!=\! c_{\RR} \!=\! 1/4$,
equations (\ref{eq:cweno.side})--(\ref{eq:cweno.exact}) yield
\bea
P_{\CC}(x) & = & 
       \uj - \frac{1}{12}(\ujp - 2\uj + \ujm) + \nonumber \\
  && + \frac{\ujp-\ujm}{2 \Delta x}(x-x_{j}) + 
       \frac{\ujp-2\uj+\ujm}{\Delta x^2}(x-x_j)^{2}.
\eea

In smooth regions, the coefficients $w_{i}$ of the convex combination in
(\ref{eq:cweno.reconstrcution}) are chosen to guarantee the maximum order
of accuracy (in this particular case -- order three), but in the
presence of a discontinuity they are automatically switched to the best
one-sided stencil (which generates the least oscillatory reconstruction).
The weights are taken as
\be
w_{i} = \frac{\alpha_{i}}{\sum\limits_{i} \alpha_{i}},  \quad
\mbox{where} \quad
\alpha_{i} = \frac{c_{i}}{(\epsilon + IS_{i})^{p}}, \qquad
i \in \{\CC,\RR,\LL\}.
\label{eq:cweno.alpha}
\ee
The constant $\epsilon$ guarantees that the denominator
does not vanish and is taken as $\epsilon = 10^{-6}$.
The value of $p$ may be chosen to provide the highest
accuracy in smooth areas and ensure the non-oscillatory nature of
the solution near the discontinuities
(consult \cite{jiang-shu:weno}, see also
\cite{levy-puppo-russo:1d, levy-puppo-russo:balance}). In
\cite{jiang-shu:weno} the value $p=2$ was empirically selected,
and here we use the same $p$ in most of the examples presented below.
Finally, the smoothness indicators, $IS_{i}$, are defined as
\[
IS_{i}=\sum_{l=1}^{2}
\int\limits_{x_{j-1/2}}^{x_{j+1/2}}
{(\Delta x)}^{2l-1} (P_{i}^{(l)}(x))^{2} dx.
\]
A direct computation then results with 
\bea
 IS_{\LL} & = &  (\uj - \ujm)^2, \qquad
 IS_{\RR} = (\ujp - \uj)^2, \nonumber \\
 IS_{\CC} & = & \frac{13}{3}(\ujp-2\uj+\ujm)^2 + \frac{1}{4}(\ujp - \ujm)^2.
  \label{eq:cweno.is}
\eea
It is easy to see that in the presence of large gradients, this reconstruction
switches to one of the second-order one-sided linear reconstructions, $P_{\RR}$ 
or $P_{\LL}$.  For more details we refer to \cite{levy-puppo-russo:balance}.


\section{The Fully-Discrete One-Dimensional Construction}
\label{section:fully-discrete}
\setcounter{equation}{0}
\setcounter{figure}{0}
\setcounter{table}{0}

In this section we present the new third-order method in the fully-discrete
framework. Since we are mainly interested in deriving the semi-discrete scheme,
we will concentrate only on the details which are required for that task.  
The scheme we derive here,
is a third-order extension of the fully-discrete second-order scheme presented
in \cite{kurganov-tadmor:semi}.

Following \cite{kurganov-tadmor:semi}, we would like to augment the 
integration over the Riemann fans by a more accurate information about the
{\em local\/} speed of wave propagation.
We start by assuming that in every cell, $I_{j}$, we have reconstructed a 
piecewise polynomial interpolant, $P_{j}(x,t^{n})$, from the previously
computed cell averages, $\{\uj\}$. Then, an upper bound on the speed of
propagation of discontinuities at the cell boundaries, $x_{j+1/2}$,
is given by
\begin{equation}
a_{j+1/2}^n=
\max_{u\in{\cal C}(u_{j+1/2}^-,u_{j+1/2}^+)}
\rho\Bigl(\frac{\partial f}{\partial u}(u)\Bigr),
\label{speed}
\end{equation}
where $\rho(A)$ denotes the spectral radius of a matrix $A$, i.e.,
 $\rho(A):=\max\limits_i|{\lambda}_i(A)|$, with ${\lambda}_i(A)$ being its eigenvalues.
We denote by $u_{j+1/2}^{+}$ and $u_{j+1/2}^{-}$ the left and right
intermediate values of $u(x,t^{n})$ at $x_{j+1/2}$, i.e.,
\[
u_{j+1/2}^{+} := P_{j+1}(x_{j+1/2},t^n), \quad
u_{j+1/2}^{-} := P_{j}(x_{j+1/2},t^n),
\]
and by ${\cal C}(u_{j+1/2}^-,u_{j+1/2}^+)$ a curve
in phase space that connects $u_{j+1/2}^-$ and $u_{j+1/2}^+$ via the Riemann fan.
\begin{remark}
In most practical applications, these local maximal speeds can
be easily evaluated. For example, in the genuinely nonlinear or linearly
degenerate case one finds that (\ref{speed}) reduces to
\begin{equation}
a_{j+1/2}^n:=
\max\left\{\rho\Bigl(\frac{\partial f}{\partial u}
(u_{j+1/2}^-)\Bigr),
\rho\Bigl(\frac{\partial f}{\partial u}
(u_{j+1/2}^+)\Bigr)\right\}.
\label{aj}
\end{equation}
\end{remark}

Given the piecewise polynomial interpolant at time $t^n$, $\{P_j(x,t^n)\},$
and the local speeds of propagation, $\{a_{j+1/2}^{n}\}$, we construct the
fully-discrete, central method in two steps, which are schematically
described in Figure \ref{figure:modified}.
First, we integrate over the control volumes, 
$[x_{j-1/2,l}^{n}, x_{j-1/2,r}^{n}] \times [t^{n},{t^{n+1}}]$,
$[x_{j-1/2,r}^{n}, x_{j+1/2,l}^{n}] \times [t^{n},{t^{n+1}}]$,
and $[x_{j+1/2,l}^{n}, x_{j+1/2,r}^{n}] \times [t^{n},{t^{n+1}}]$,
obtaining $\bar{w}_{j-1/2}^{n+1}$, $\bar{w}_{j}^{n+1}$ and
$\bar{w}_{j+1/2}^{n+1}$, respectively.  Due to the
finite speed of propagation, the points $x_{j+1/2,l}^n$ and $x_{j+1/2,r}^n$,
\[
x_{j+1/2,l}^n:=x_{j+1/2}-a_{j+1/2}^n\Delta t, \quad
x_{j+1/2,r}^n:=x_{j+1/2}+a_{j+1/2}^n\Delta t,
\]
separate between smooth and non-smooth regions. That is, the solution of
equation (\ref{CL}) subject to the piecewise polynomial initial data
prescribed at time $t=t^n$,
may be non-smooth only inside the intervals
$[x_{j+1/2,l}^n\,,\,x_{j+1/2,r}^n]$ for $t\in[t^n,t^{n+1})$.

In the second step,
we repeat the non-oscillatory reconstruction (this time on 
a nonuniformly spaced grid) and project the obtained reconstruction
on the original, uniform grid, ending up with the cell averages
at the next time level $t^{n+1}$, $\{\bar{u}_{j}^{n+1}\}$.
This last step does not involve
time integration, and was introduced in the context of changing staggered
methods into non-staggered methods in \cite{jiang-levy-osher-tadmor:stg}.

\begin{figure}[H]
\begin{center}
\hspace*{-0.5cm}
Contact Author for Figure If Necessary
\caption{Modified Central Differencing}
\label{figure:modified}
\end{center}
\end{figure}

We now turn to the detailed description of this algorithm.
Assume that the piecewise polynomial reconstruction in cell $I_{j}$
at time $t^n$ is of the form
\be
P_{j}(x,t^{n}) = A_{j} + B_{j}(x-x_{j}) + \frac{1}{2}C_{j}(x-x_{j})^{2}.
\label{eq:reconstruction.p}
\ee
Then a direct computation of the integrals over the control volumes,
$[x_{j+1/2,l}^{n}, x_{j+1/2,r}^{n}] \times [t^{n},{t^{n+1}}]$ and
$[x_{j-1/2,r}^{n}, x_{j+1/2,l}^{n}] \times [t^{n},{t^{n+1}}]$, yields
\bea
\lefteqn{\bar{w}_{j+1/2}^{n+1} = \frac{A_{j}+A_{j+1}}{2} 
+\frac{\Delta x-a_{j+1/2}^n\Delta t}{4}(B_{j}-B_{j+1})+}
\label{wjhalf} \\
&&+\left(\frac{\Delta x^{2}}{16} - 
\frac{a_{j+1/2}^n \Delta t \Delta x}{8} + 
\frac{(a_{j+1/2}^n\Delta t)^{2}}{12}\right)(C_{j}+C_{j+1}) - \nonumber \\
&&-\frac{1}{2 a_{j+1/2}^n\Delta t}
\left\{\int\limits_{t^{n}}^{t^{n+1}} \left[ f(u(x_{j+1/2,r}^n,t))dt - 
f(u(x_{j+1/2,l}^n,t)) \right]dt \right\};\nonumber
\eea
and
\bea
\lefteqn{\bar{w}_{j}^{n+1}=A_{j}+\frac{\Delta t}{2}
(a_{j-1/2}^n-a_{j+1/2}^n)B_{j}+}\label{wj} \\
&&\!\!\!\!\!\!\!\!\!\!\!\!\!\!\!\!\!\!\!\!\!\!\!\!\!\!\!
+\left[\frac{(\Delta x)^{2}}{24}-\frac{\Delta t \Delta x}{12}
(a_{j-1/2}^n+a_{j+1/2}^n)
+\frac{(\Delta t)^{2}}{6}\Bigl({(a_{j-1/2}^n)}^{2}-
a_{j-1/2}^na_{j+1/2}^n+{(a_{j+1/2}^n)}^{2}\Bigr)\right]C_{j}-\nonumber \\
&&\!\!\!\!\!\!-\frac{1}{\Delta x-\Delta t(a_{j-1/2}^n-a_{j+1/2}^n)}
\left\{\int\limits_{t^{n}}^{t^{n+1}}\left[ f(u(x_{j+1/2,l}^n,t))dt-
f(u(x_{j-1/2,r}^n,t)) \right]dt \right\},\nonumber
\eea
respectively. To complete these computations, one should approximate
the flux integrals on the RHS of (\ref{wjhalf}) and (\ref{wj}) using, e.g.,
Simpson's quadrature as described in \S\ref{section:central}.

At this stage, the approximate cell averages, 
 $\{\bar{w}_{j+\frac{1}{2}}^{n+1},\,\bar{w}_j^{n+1}\},$ 
realize the solution at $t=t^{n+1}$ over a nonuniform grid,
which is oversampled by twice the number of the original cells at $t=t^n$.
To convert these nonuniform averages back into the original grid,
we proceed along the lines of \cite{jiang-levy-osher-tadmor:stg}.

First, from the cell averages,
$\bar{w}_{j+\frac{1}{2}}^{n+1},\,\bar{w}_j^{n+1}$,
given by (\ref{wjhalf})--(\ref{wj}), we reconstruct a third-order,
piecewise polynomial,
non-oscillatory interpolant (e.g., the CWENO interpolant described
in \S\ref{subsection:cweno}), which we will denote by
$\tilde{w}_{j+1/2}^{n+1}(x)$ and $\tilde{w}_{j}^{n+1}(x)$, respectively.
In fact, we do not need any high-order reconstruction
$\tilde{w}_{j}^{n+1}(x)$ since it will be averaged out (consult
Figure \ref{figure:modified}).

We note in passing that even for a nonuniform grid data, the CWENO
interpolant can be written explicitly
(in the spirit of \S\ref{subsection:cweno}), but these details are
irrelevant for the semi-discrete scheme, which will be described in
\S\ref{section:semi-discrete}. At that point, all that we need is
to assume that such a reconstruction exists and that
for all $j$ it takes the form
\be
\begin{array}{lll}
\tilde{w}_{j+1/2}^{n+1}(x)&=&\displaystyle{\tilde{A}_{j+1/2}+
\tilde{B}_{j+1/2}(x-x_{j+1/2})+\frac{1}{2}\tilde{C}_{j+1/2}(x-x_{j+1/2})^{2}},
\\ \\
\tilde{w}_{j}^{n+1}(x)&=&{\bar w}_{j}^{n+1},
\end{array}
\label{eq:reconstruction.w}
\ee
in the non-smooth region, $(x_{j+1/2,l}^n,x_{j+1/2,r}^n),$ and in the smooth region,
$(x_{j-1/2,r}^n,x_{j+1/2,l}^n),$ respectively.
Given (\ref{wjhalf}), (\ref{wj}) and (\ref{eq:reconstruction.w}),
we conclude by computing the new cell averages at 
time $t^{n+1}$ according to
\bea
\bar{u}_{j}^{n+1}&=&\frac{1}{\Delta x}\left[
\int\limits_{x_{j-1/2}}^{x_{j-1/2,r}^n}\tilde{w}_{j-1/2}^{n+1}(x)dx+
\int\limits_{x_{j-1/2,r}^n}^{x_{j+1/2,l}^n}\tilde{w}_{j}^{n+1}(x)dx+
\int\limits_{x_{j+1/2,l}^n}^{x_{j+1/2}}\tilde{w}_{j+1/2}^{n+1}(x)dx\right]=
\nonumber\\
&=&\lambda a_{j-1/2}^n\tilde{A}_{j-1/2}+
\left[1-\lambda(a_{j-1/2}^n+a_{j+1/2}^n)\right]{\bar w}_{j}^{n+1}
+\lambda a_{j+1/2}^n\tilde{A}_{j+1/2}+\nonumber\\
&&+\frac{\lambda\Delta t}{2}\left({(a_{j-1/2}^n)}^{2}\tilde{B}_{j-1/2}-
{(a_{j+1/2}^n)}^{2}\tilde{B}_{j+1/2}\right)+\nonumber\\
&&+\frac{\lambda(\Delta t)^{2}}{6}\left({(a_{j-1/2}^n)}^{3}\tilde{C}_{j-1/2}+
{(a_{j+1/2}^n)}^{3}\tilde{C}_{j+1/2}\right).
\label{eq:fully.scheme}
\eea

\begin{remark}
The third-order reconstruction (\ref{eq:reconstruction.w}) is necessary
in order to guarantee the overall third-order accuracy, since simple averaging
over $[x_{j-\frac{1}{2}},x_{j+\frac{1}{2}}]$ (without reconstruction)
reduces the order of the resulting scheme
(see \cite{jiang-levy-osher-tadmor:stg}).
\end{remark}


\section{The Semi-Discrete Scheme}         \label{section:semi-discrete}
\setcounter{equation}{0}
\setcounter{figure}{0}
\setcounter{table}{0}

We are now ready to derive our main result, which is the new third-order,
semi-discrete, central scheme.  First, we describe our ideas in the one-dimensional
framework and then we extend them to multidimensional problems.

\subsection{One-Dimensional Problems} \label{subsec:1D}

We start with the derivation of the third-order semi-discrete scheme
for one-dimensional (systems of) hyperbolic conservation laws.
Using the fully-discrete scheme
obtained in \S\ref{section:fully-discrete}, the semi-discrete approximation
can be directly written as the limit
\be
\frac{d}{dt}{\bar u}_{j}(t)=
\lim_{\Delta t\rightarrow 0} \frac{{\bar u}_j^{n+1}-\uj}{\Delta t}.
\label{eq:semi.start}
\ee
Substituting (\ref{eq:fully.scheme}) into (\ref{eq:semi.start}) results with
\bea
\frac{d{\bar u}_{j}}{dt}&=&\lim_{\Delta t\rightarrow 0}\left\{
\frac{1}{\Delta x}a_{j-1/2}^n\tilde{A}_{j-1/2} 
-\frac{1}{\Delta x}(a_{j-1/2}^n+a_{j+1/2}^n){\bar w}_{j}^{n+1}+
\frac{1}{\Delta x}a_{j+1/2}^n\tilde{A}_{j+1/2}+\right.\nonumber\\
&&~~~~~\left.+\frac{1}{\Delta t}({\bar w}_{j}^{n+1}-\uj)\right\}.
\label{eq:dot.u.2}
\eea
In the limit as $\Delta t \rightarrow 0$, all the Riemann fans have
zero widths and therefore,
\be
\tilde{A}_{j+1/2} = {\bar w}_{j+1/2}^{n+1}, \qquad
\tilde{A}_{j-1/2} = {\bar w}_{j-1/2}^{n+1}.
\label{a=w}
\ee
Using (\ref{eq:reconstruction.p}) we can also obtain
\bean
u(x_{j+1/2,r}^n,t) &\longrightarrow& P_{j+1}(x_{j+1/2},t)=
A_{j+1}-\frac{\Delta x}{2}B_{j+1}+\frac{{(\Delta x)}^2}{8}C_{j+1}=
:u_{j+1/2}^{+}(t),\\
u(x_{j+1/2,l}^n,t) &\longrightarrow& P_{j}(x_{j+1/2},t)=
A_j+\frac{\Delta x}{2}B_j+\frac{{(\Delta x)}^2}{8}C_j=
:u_{j+1/2}^{-}(t).
\eean
Finally, plugging (\ref{wjhalf}), (\ref{wj}) and (\ref{a=w})
into (\ref{eq:dot.u.2}) we compute the time limit
explicitly, ending up with our new semi-discrete scheme,
\be
\begin{array}{lll}
\displaystyle{\frac{d{\bar u}_{j}}{dt}}&=&
\displaystyle{-\frac{1}{2\Delta x}
\left[f(u_{j+1/2}^{+}(t))+f(u_{j+1/2}^{-}(t))-
f(u_{j-1/2}^{+}(t))-f(u_{j-1/2}^{-}(t))\right]+}\\ \\
&&\displaystyle{+\frac{a_{j+1/2}(t)}{2\Delta x}
\left[u_{j+1/2}^{+}(t)-u_{j+1/2}^{-}(t)\right]-
\frac{a_{j-1/2}(t)}{2\Delta x}\left[u_{j-1/2}^{+}(t)-u_{j-1/2}^{-}(t)\right]},
\end{array}
\label{semidiscrete_scheme}
\ee
with local speeds $a_{j+1/2}(t)$, e.g., $a_{j+1/2}(t):=
\max\left\{\rho\Bigl(\frac{\partial f}{\partial u}
(u_{j+1/2}^-(t))\Bigr),
\rho\Bigl(\frac{\partial f}{\partial u}
(u_{j+1/2}^+(t))\Bigr)\right\}.$

\begin{remarks}
\begin{enumerate}
\item Our third-order scheme, (\ref{semidiscrete_scheme}),
admits the conservative form,
\be
\frac{d{\bar u}_{j}}{dt}=-\frac{H_{j+1/2}(t)-H_{j-1/2}(t)}{\Delta x},
\label{semi 2nd_order flux_form}
\ee
with the numerical flux
\be
H_{j+1/2}(t):=
\frac{f(u_{j+1/2}^+(t))+f(u_{j+1/2}^-(t))}{2}-
\frac{a_{j+1/2}(t)}{2}\left[u_{j+1/2}^+(t)-
u_{j+1/2}^-(t)\right].
\label{flux}
\ee
This scheme is a natural generalization of the
second-order semi-discrete scheme from \cite{kurganov-tadmor:semi}.
Moreover, the second-order scheme has exactly the same form,
(\ref{semi 2nd_order flux_form})--(\ref{flux});
the only difference is in the more accurate
computation of the intermediate values, $u_{j+1/2}^+(t)$ and
$u_{j+1/2}^-(t)$. It is interesting to note that also the
fully-discrete, staggered, second- and third-order central schemes have
the same structure (see \cite{liu-tadmor:3rd}).
\item Similar to the case of the second-order scheme,
\cite{kurganov-tadmor:semi}, the non-oscillatory property of
the piecewise parabolic reconstruction,
(\ref{eq:reconstruction.p}), will guarantee
the non-oscillatory nature of our semi-discrete scheme. But unlike
the piecewise linear reconstruction utilized in the second-order
method, a piecewise parabolic reconstruction can be only
essentially non-oscillatory.
This means that, in principle, such a reconstruction may increase
the total variation of the computed piecewise constant solution.
Our numerical examples, however, demonstrate that the growth of the
total variation is always bounded.
Such desirable behavior of bounded total variation
in the context of central-WENO schemes, was
already observed in \cite{levy-puppo-russo:tv}.
\item We would like to stress once again the simplicity of our new
method, which does not require any (approximate) Riemann solver
or any use of the characteristic variables -- the reconstruction of
piecewise polynomial interpolant, (\ref{eq:reconstruction.p}), is
carried out {\em componentwise}.  In particular, unlike the standard
central schemes, but similar to the second-order semi-discrete method
in \cite{kurganov-tadmor:semi}, our method is based on one grid
(and not on staggering between two grids).
This can be a big advantage (compared with the traditional central schemes)
when dealing with boundary conditions and complex geometries.
\end{enumerate}
\end{remarks}

\bigskip
Next, let us consider the general convection-diffusion equation, (\ref{CD}).
Similar to the case of the second-order semi-discrete scheme,
\cite{kurganov-tadmor:semi},
operator splitting is not needed. We can apply our third-order
semi-discrete scheme, (\ref{semi 2nd_order flux_form})--(\ref{flux}),
to the (degenerate) parabolic equation, (\ref{CD}),
in a straightforward manner. This results in the following scheme,
\begin{equation}
\frac{d{\bar u}_{j}}{dt}=-\frac{H_{j+1/2}(t)-H_{j-1/2}(t)}{\Delta x}
+Q_j(t).
\label{semi_CD}
\end{equation}
Here, $H_{j+1/2}(t)$ is our numerical convection flux, (\ref{flux}),
and $Q_j(t)$ is a high-order approximation to the
diffusion term, ${Q(u,u_x)}_x$. In the examples below we use the fourth-order
central differencing of the form
\bea
Q_j(t)=\frac{1}{12\Delta x}
\Bigl[\!\!\!\!&-&\!\!\!\!Q(u_{j+2}(t),{(u_x)}_{j+2,j})+
8Q(u_{j+1}(t),{(u_x)}_{j+1,j})-\nonumber\\
&&-8Q(u_{j-1}(t),{(u_x)}_{j-1,j})+Q(u_{j-2}(t),{(u_x)}_{j-2,j})\Bigr],
\label{diff_flux}
\eea
where
\bea
{(u_x)}_{j+2,j}\!\!\!&:=&\!\!\!\frac{1}{12\Delta x}
\Bigl[25u_{j+2}(t)-48u_{j+1}(t)+36u_j-16u_{j-1}(t)+3u_{j-2}(t)\Bigr],
\nonumber \\
{(u_x)}_{j+1,j}\!\!\!&:=&\!\!\!\frac{1}{12\Delta x}
\Bigl[3u_{j+2}(t)+10u_{j+1}(t)-18u_j+6u_{j-1}(t)-u_{j-2}(t)\Bigr],\nonumber \\
{(u_x)}_{j-1,j}\!\!\!&:=&\!\!\!\frac{1}{12\Delta x}
\Bigl[u_{j+2}(t)-6u_{j+1}(t)+18u_j-10u_{j-1}(t)-3u_{j-2}(t)\Bigr],\nonumber \\
{(u_x)}_{j-2,j}\!\!\!&:=&\!\!\!\frac{1}{12\Delta x}
\Bigl[-3u_{j+2}(t)+16u_{j+1}(t)-36u_j+48u_{j-1}(t)-25u_{j-2}(t)\Bigr];
\label{u_x}
\eea
and $\{u_j(t)\}$ are point-values of the reconstructed polynomials,
(\ref{eq:reconstruction.p}), i.e., $u_j(t)=P_j(x_j,t)$.

\subsection{Multi-Dimensional Extensions} \label{subsec:2D}

Without loss of generality, let us consider the
{\em two-dimensional} (system of) convection-diffusion equations,
\begin{equation}
u_t+f(u)_x+g(u)_y={Q^x(u,u_x,u_y)}_x+{Q^y(u,u_x,u_y)}_y,
\label{2D CL_CD}
\end{equation}
where the case $Q^x\equiv Q^y\equiv 0$ corresponds to the 2D pure hyperbolic
problem.

Suppose that we have computed an approximate solution to
(\ref{2D CL_CD}) at some time $t$, and have reconstructed
a two-dimensional piecewise polynomial, third-order,
essentially non-oscillatory interpolant over the uniform spatial grid,
$(x_j,y_k)=(j\Delta x,k\Delta y)$.

Following \cite{kurganov-tadmor:semi}, 
the 2D extension of our third-order semi-discrete scheme,
(\ref{semi_CD}),(\ref{flux}), can be written in the following form,
\bea
\frac{d{\bar u}_{j,k}}{dt}& =& -\frac{H_{j+1/2,k}^x(t)-H_{j-1/2,k}^x(t)}{\Delta x}
-\frac{H_{j,k+1/2}^y(t)-H_{j,k-1/2}^y(t)}{\Delta y} + \label{semi multi_CD} \\
&&+Q_{j,k}^x(t)+Q_{j,k}^y(t). \nonumber
\eea
Here, $H^x_{j+1/2,k}(t)$ and $H^y_{j,k+1/2}(t)$ are
$x$- and $y$-numerical convection fluxes, respectively (they can be viewed
as a generalization of the one-dimensional flux, (\ref{flux})),
\begin{eqnarray}
H^x_{j+1/2,k}(t)\!\!\!&:=&\!\!\!
\frac{f(u_{j+1/2,k}^+(t))+f(u_{j+1/2,k}^-(t))}{2}- \nonumber \\
&&-\frac{a^x_{j+1/2,k}(t)}{2}\left[u_{j+1/2,k}^+(t)-
u_{j+1/2,k}^-(t)\right],\nonumber\\
\label{2D_fluxes}\\
H^y_{j,k+1/2}(t)\!\!\!&:=&\!\!\!
\frac{g(u_{j,k+1/2}^+(t))+g(u_{j,k+1/2}^-(t))}{2}- \nonumber \\
&&-\frac{a^y_{j,k+1/2}(t)}{2}\left[u_{j,k+1/2}^+(t)-
u_{j,k+1/2}^-(t)\right].\nonumber
\end{eqnarray}
The numerical fluxes, (\ref{2D_fluxes}), are expressed in terms of the 
intermediate values, $u_{j+1/2,k}^{\pm}(t),\\u_{j,k+1/2}^{\pm}(t)$,
which are obtained from the piecewise
polynomial reconstruction. The local speeds, $a^x_{j+1/2,k}(t)$ and
$a^y_{j,k+1/2}(t),$ are computed, e.g., by
\begin{eqnarray}
a^x_{j+1/2,k}(t)\!\!&:=&\!\!
\max\left\{\rho\Bigl(\frac{\partial f}{\partial u}
(u_{j+1/2,k}^-(t))\Bigr),
\rho\Bigl(\frac{\partial f}{\partial u}
(u_{j+1/2,k}^+(t))\Bigr)\right\},\nonumber\\
\label{2D_a}\\
a^y_{j,k+1/2}(t)\!\!&:=&\!\!
\max\left\{\rho\Bigl(\frac{\partial g}{\partial u}
(u_{j,k+1/2}^-(t))\Bigr),
\rho\Bigl(\frac{\partial g}{\partial u}
(u_{j,k+1/2}^+(t))\Bigr)\right\}.\nonumber
\end{eqnarray}
Finally, $Q_{j,k}^x(t)$  and $Q_{j,k}^y(t)$ are high-order,
central differencing approximations to the diffusion terms,
${Q^x(u,u_x,u_y)}_x$ and ${Q^y(u,u_x,u_y)}_y$.

\begin{remarks}
\begin{enumerate}
\item We would like to emphasize that the problem of constructing a two-dimensional,
third-order, non-oscillatory interpolant is highly non-trivial.
Several essentially 2D reconstructions were proposed in
\cite{levy-puppo-russo:2d-3,
levy-puppo-russo:balance, levy-puppo-russo:2d-4}.
Alternatively, one can use one-dimensional CWENO reconstruction, 
direction by direction, in order to
compute the intermediate values,
$u_{j+1/2,k}^{\pm}(t)$ and $u_{j,k+1/2}^{\pm}(t)$.

Following is the recipe for the computation of $u_{j+1/2,k}^+$
(the computation of other intermediate values can be carried out in
the similar way).
\be
u_{j+1/2,k}^+=w_{\LL}P_{\LL}^k(x_{j+1/2})+w_{\RR}P_{\RR}^k(x_{j+1/2})
+w_{\CC}P_{\CC}^k(x_{j+1/2}),
\label{dim*dim_1}
\ee
where the $P$'s are the polynomials introduced in \S \ref{subsection:cweno},
\bea
P_{\RR}^k(x)&=&{\bar u}_{j,k}+
\frac{{\bar u}_{j+1,k}-{\bar u}_{j,k}}{\Delta x}(x-x_{j}),
\quad
P_{\LL}^k(x)={\bar u}_{j,k}+
\frac{{\bar u}_{j,k}-{\bar u}_{j-1,k}}{\Delta x}(x-x_{j}),\nonumber\\
P_{\CC}^k(x)&=&
{\bar u}_{j,k}-\frac{1}{12}({\bar u}_{j+1,k}-2{\bar u}_{j,k}+{\bar u}_{j-1,k})-
\frac{1}{12}({\bar u}_{j,k+1}-2{\bar u}_{j,k}+{\bar u}_{j,k-1})+
\nonumber \\
&&+\frac{{\bar u}_{j+1,k}-{\bar u}_{j-1,k}}{2 \Delta x}(x-x_{j}) +
\frac{{\bar u}_{j+1,k}-2{\bar u}_{j,k}+{\bar u}_{j-1,k}}{\Delta x^2}(x-x_j)^{2}.
 \label{dim*dim_2} 
\eea
The weights, $w_{\LL}, w_{\RR}, w_{\CC}$ which are given by (\ref{eq:cweno.alpha}),
and are based on the smoothness indicators in (\ref{eq:cweno.is}).

Note that the only difference between this reconstruction and the
1D reconstruction, (\ref{eq:cweno.reconstrcution})--(\ref{eq:cweno.is}),
is an additional term in $P_{\CC}^k(x)$,
$-\frac{1}{12}({\bar u}_{j,k+1}-2{\bar u}_{j,k}+{\bar u}_{j,k-1})$,
which corresponds to the second derivative
in the $y$ direction and guarantees
the third-order accuracy of the computed intermediate values.

This `dimension by dimension' approach, was implemented in Example 5 below.

\item It is straightforward to extend the two-dimensional scheme, (\ref{semi multi_CD}),
to more space dimensions.  In particular, the dimension-by-dimension approach is
a very simple and promising approach for multi-dimensional problems.
\end{enumerate}
\end{remarks}


\section{Numerical Examples}           \label{section:numerical}
\setcounter{equation}{0}
\setcounter{figure}{0}
\setcounter{table}{0}

We conclude the paper with a number of numerical examples. 
Here, in order to retain the overall high accuracy,
the semi-discrete scheme is combined
with a high-order, stable ODE solver to complete the spatio-temporal
discretization.  Numerically, we observed that a variety of
explicit methods provide satisfactory results in the context of our
semi-discrete scheme.

For the inviscid problems (Examples 1, 2, 3 and 5), we used
the third-order total variation diminishing (TVD) Runge-Kutta type method 
introduced by Shu and Osher in \cite{SO}.
However, if we apply this time-integration method or any other
standard Runge-Kutta type method to (degenerate) parabolic
problems, the time-step can be very small due to their strict 
stability restrictions.

To overcome this difficulty, we used (in Examples 4 and 5)
the third-order ODE solver (called DUMKA3) by Medovikov, \cite{Med}. 
This explicit method has larger stability
domains (compared with the standard Runge-Kutta methods), which allow
larger time-steps. In practice, DUMKA3 works as fast as
implicit methods (see \cite{Med} for details).

We abbreviate by SD3 our third-order semi-discrete scheme,
which will be combined with the third-order TVD Runge-Kutta type method
(RK3) or with DUMKA3.


\subsection*{Example 1: Linear Accuracy Test}

Consider the scalar linear hyperbolic equation
\be
u_t+u_x=0, \qquad x\in[0,2\pi],
\label{lin eq}
\ee
augmented with the smooth initial data, $u(x,0)=\sin x$, and periodic boundary
conditions.
This simple problem admits a global classical solution, which was
computed at time $T=1$ with a varying number of grid points, $N$.

In Table \ref{table:linear} we check
the accuracy of our third-order semi-discrete scheme,
SD3, coupled with the RK3 ODE solver. If instead of computing
the approximate convergence rate between two consecutive mesh refinings,
one approximates the convergence rate between $N=40$ and
$N=1280$, the results are 3.27 in the $L^{1}$-norm and 2.91 in the
$L^{\infty}$-norm. This clearly demonstrates that our scheme is third-order.

The error is measured in terms of the pointwise values,
\[
{\|{\tilde u}-u\|}_{L^1}:=\Delta x \sum\limits_j |{\tilde u}_j(T)-u(x_j,T)|,
\quad
{\|{\tilde u}-u\|}_{L^{\infty}}:=\max_j |{\tilde u}_j(T)-u(x_j,T)|.
\]
Here, ${\tilde u}$ is an approximate solution, which is realized
by its values at the grid points, $x_j$,
\[
{\tilde u}_j(T)=P_j(x_j,T),
\]
where the $P_j$'s are the piecewise parabolic interpolants,
(\ref{eq:reconstruction.p}), constructed at the final time $t=T$.

\begin{table}[!h]
\begin{center}
\begin{tabular}{|c|cccc|} \hline
&&&&\\
N & $L^1$-error & rate & $L^{\infty}$-error & rate \\
&&&&\\ \hline
&&&&\\
40 & 4.492e-02 & -- & 2.822e-02 & -- \\
80 & 1.092e-02 & 2.04 & 1.065e-02 & 1.41 \\
160 & 2.162e-03 & 2.34 & 3.426e-03 & 1.64 \\
320 & 1.811e-04 & 3.58 & 4.705e-04 & 2.86 \\
640 & 9.267e-06 & 4.29 & 2.267e-05 & 4.38 \\
1280 & 5.409e-07 & 4.10 & 1.171e-06 & 4.27 \\
&&&&\\ \hline
\end{tabular}
\caption{Accuracy test for the linear advection problem, (\ref{lin eq}); 
 The errors at $T=1$.
\label{table:linear}}
\end{center}
\end{table}


\subsection*{Example 2: Burgers' Equation}

In this example we approximate solutions to the inviscid
Burgers' equation,
\begin{equation}
u_t+{\left(\frac{u^2}{2}\right)}_x=0, \qquad x\in[0,2\pi],
\label{Burgers}
\end{equation}
augmented with the smooth initial data, $u(x,0)=0.5+\sin x$,
and periodic boundary conditions.

The unique entropy solution of (\ref{Burgers}) develops a
shock discontinuity at the critical time $T_c=1$.
Table \ref{table:burgers} shows the
$L^1$- and $L^{\infty}$-norms of the errors at the pre-shock time
$T=0.5$, when the solution is still smooth. 
Once again, when approximating the convergence rate by looking at
the errors for $N=40$ and $N=1280$, we get 3.25 in the $L^{1}$-norm
and 3.10 in the $L^{\infty}$-norm.
These results indicate that our method is also third-order accurate 
when the accuracy is measured in nonlinear problems.

In Figures \ref{SD3_Burg_40}--\ref{SD3_Burg_80} we present the approximate
solutions at the post-shock
time $T=2$, when the shock is well developed.
The essentially non-oscillatory nature of our scheme can
be clearly observed.

\begin{table}[!h]
\begin{center}
\begin{tabular}{|c|cccc|} \hline
&&&&\\
N & $L^1$-error & rate & $L^{\infty}$-error & rate \\
&&&&\\ \hline
&&&&\\
40 & 2.370e-02 & -- & 2.225e-02 & -- \\
80 & 5.759e-03 & 2.04 & 9.053e-03 & 1.30 \\
160 & 1.161e-03 & 2.31 & 2.921e-03 & 1.63 \\
320 & 9.541e-05 & 3.61 & 3.926e-04 & 2.90 \\
640 & 4.882e-06 & 4.29 & 1.778e-05 & 4.46 \\
1280 & 3.044e-07 & 4.00 & 5.732e-07 & 4.96 \\
&&&&\\ \hline
\end{tabular}
\caption{Accuracy test for Burgers equation, (\ref{Burgers}); The pre-shock errors \label{table:burgers}}
\end{center}
\end{table}

\begin{center}
\begin{minipage}[t]{7.5cm}
\begin{figure}[H]
Contact Author for Figure If Necessary
\caption{Burgers equation, (\ref{Burgers}); $T=2$, {\bf N=40}.
\label{SD3_Burg_40}}
\end{figure}
\end{minipage}\ \hspace{0.4cm}\
\begin{minipage}[t]{7.5cm}
\begin{figure}[H]
Contact Author for Figure If Necessary
\caption{Burgers equation, (\ref{Burgers}); $T=2$, {\bf N=80}.
\label{SD3_Burg_80}}
\end{figure}
\end{minipage}
\end{center}

\subsection*{Example 3: Euler Equations of Gas Dynamics}
 
Let us consider the one-dimensional Euler system,
\[
\frac{\partial}{\partial t} \left[ \begin{array}{l} \rho\\ m\\ E \end{array}
\right]+
\frac{\partial}{\partial x} \left[ \begin{array}{c} m\\ \rho u^2+p\\ u(E+p)
\end{array} \right]=0,
\qquad p=(\gamma-1) \cdot \left(E-\frac{\rho}{2}u^2\right),
\]
where $\,\rho,\, u,\, m=\rho u,\, p\,$ and $\,E\,$ are the density, velocity,
momentum, pressure and total energy, respectively. 
We solve this system subject to Sod's Riemann initial data, 
proposed in \cite{sod:hyperbolic},
\[
\vec u(x,0)=
\left\{
\begin{array}{ll}
\vec u_L={(1,0,2.5)}^T, &~~x<0,\\
\vec u_R={(0.125,0,0.25)}^T, &~~x>0.
\end{array}
\right .
\]

The approximations to the density, velocity and pressure obtained by the
SD3 scheme with the RK3 time discretization are presented
in Figures \ref{Sod_den_200}--\ref{Sod_pre_400}.
The coefficient $p$ in the smoothness indicator,
(\ref{eq:cweno.alpha})--(\ref{eq:cweno.is}),
was taken as $0.6$, which seems to be the optimal value in this
specific example.

We would like to stress again that our SD3 scheme does not
require the characteristic decomposition. To improve the
resolution of the contact discontinuity, which is always
smeared while the solution to the system of Euler equations is
computed by the central method, we implemented the Artificial Compression
Method (ACM) by Harten, \cite{harten:ACM}.
In the context of central schemes, the ACM can be implemented
as a corrector step to the component-wise approach
(see \cite{nessyahu-tadmor:non-oscillatory} for details).

\begin{center}
\begin{minipage}[t]{7.5cm}
\begin{figure}[H]
\hspace*{-0.5cm}
Contact Author for Figure If Necessary
\caption{~Sod problem -- {\em density}. {\bf N=200, T=0.1644}.}
\label{Sod_den_200}
\end{figure}
\end{minipage}\ \hspace{0.6cm}\
\begin{minipage}[t]{7.5cm}
\begin{figure}[H]
\hspace*{-0.5cm}
Contact Author for Figure If Necessary
\caption{Sod problem -- {\em density}. {\bf N=400, T=0.1644}.}
\label{Sod_den_400}
\end{figure}
\end{minipage}
\end{center}

\begin{center}
\begin{minipage}[t]{7.5cm}
\begin{figure}[H]
\hspace*{-0.5cm}
Contact Author for Figure If Necessary
\caption{Sod problem -- {\em velocity}. {\bf N=200, T=0.1644}.}
\label{Sod_vel_200}
\end{figure}
\end{minipage}\ \hspace{0.6cm}\
\begin{minipage}[t]{7.5cm}
\begin{figure}[H]
\hspace*{-0.5cm}
Contact Author for Figure If Necessary
\caption{Sod problem -- {\em velocity}. {\bf N=400, T=0.1644}.}
\label{Sod_vel_400}
\end{figure}
\end{minipage}
\end{center}
 
\begin{center}
\begin{minipage}[t]{7.5cm}
\begin{figure}[H]
\hspace*{-0.5cm}
Contact Author for Figure If Necessary
\caption{Sod problem -- {\em pressure}. {\bf N=200, T=0.1644}.}
\label{Sod_pre_200}
\end{figure}
\end{minipage}\ \hspace{0.6cm}\
\begin{minipage}[t]{7.5cm}
\begin{figure}[H]
\hspace*{-0.5cm}
Contact Author for Figure If Necessary
\caption{Sod problem -- {\em pressure}. {\bf N=400, T=0.1644}.}
\label{Sod_pre_400}
\end{figure}
\end{minipage}
\end{center}


\subsection*{Example 4: Convection-Diffusion Equations --
the Buckley-Leverett Model}
 
In this example we solve the one-dimensional Buckley-Leverett equation,
\begin{equation}
u_t+f(u)_x=\varepsilon{(\nu(u)u_x)}_x, \quad\varepsilon \nu(u)\geq 0,
\label{BL}
\end{equation}
which can be viewed as a prototype model
for the two-phase flow in oil reservoirs.
Typically, $\nu(u)$ vanishes at some values of $u$, and thus (\ref{BL})
is a degenerate parabolic equation.  Specifically, we take
\[
f(u)=\frac{u^2}{u^2+(1-u)^2}, \quad \nu(u)=4u(1-u), \quad
\varepsilon=0.01,
\]
and consider the initial value problem with the Riemann initial data,
\be
u(x,0)=\left\{
\begin{array}{lc}
0, &~~0\leq x<1-\frac{1}{\sqrt{2}},\\ \\
1, &~~1-\frac{1}{\sqrt{2}}\leq x\leq 1.
\end{array}
\right .
\label{BL_ID}
\ee
The numerical solution to this problem, obtained by the SD3 scheme
augmented with the DUMKA3 ODE solver,
is presented in Figure \ref{BL_100_800}.

\begin{center}
\begin{minipage}[t]{7.5cm}
\begin{figure}[H]
\hspace*{-0.5cm}
Contact Author for Figure If Necessary
\caption{Buckley-Leverett model, (\ref{BL})--(\ref{BL_ID}).
{\bf T=0.2}.}
\label{BL_100_800}
\end{figure}
\end{minipage}\ \hspace{0.6cm}\
\begin{minipage}[t]{7.5cm}
\begin{figure}[H]
\hspace*{-0.5cm}
Contact Author for Figure If Necessary
\caption{Buckley-Leverett model, (\ref{BL})--(\ref{BL_ID}),
including the gravitational effect, (\ref{eq:gravitation}). {\bf T=0.2}.}
\label{BL_grav_100_800}
\end{figure}
\end{minipage}
\end{center}

The model, (\ref{BL}), becomes more complicated
by adding the effects of gravitation.  This can be obtained,
e.g., by taking
\be
f(u)=\frac{u^2}{u^2+(1-u)^2}(1-5(1-u)^2),  \label{eq:gravitation}
\ee
which is non-monotone on the interval $u \in [0,1].$

The numerical solution to this initial value problem is shown
in Figure \ref{BL_grav_100_800}. Note that the exact solution to
problem  (\ref{BL})--(\ref{BL_ID}) is not available, but our
solutions seem to converge to the physically relevant solutions
in the both cases -- with gravitation or without it.

 
\subsection*{Example 5: Incompressible Euler and Navier-Stokes equations}
 
In this example we consider two-dimensional
viscous and inviscid incompressible flow governed by
the Navier-Stokes ($\nu>0$) and by the Euler ($\nu=0$) equations,
\be
{\vec u}_t+(\vec u \cdot \nabla) \vec u+\nabla p=\nu\Delta \vec u.
\label{inc_Euler}
\ee
Here, $p$ denotes the pressure, and $\vec u=(u,v)$ is the two-component
divergence-free velocity field, satisfying
\be
u_x+v_y=0.
\label{div_free}
\ee
In the 2D case (\ref{inc_Euler}) admits an equivalent scalar formulation
in terms of the vorticity,
\be
{\omega}_t+{(u\omega)}_x+{(v\omega)}_y=\nu\Delta \omega,
\label{vorticity}
\ee
where $\omega:=v_x-u_y$. The incompressibility, (\ref{div_free}), implies
that equation (\ref{vorticity}) can be written in an equivalent conservative
form,
\be
{\omega}_t+{f(\omega)}_x+{g(\omega)}_y=\nu\Delta \omega,
\label{vorticity_CL}
\ee
with a {\em global\/} convection flux, $(f,g):=(u\omega, v\omega)$.
A second-order, fully-discrete, staggered, central scheme
was used to solve the two-dimensional vorticity
equations in \cite{levy-tadmor}.  This scheme was proved to satisfy
a maximum principle for the vorticity.  (For an equivalent scheme
in the velocity formulation, see \cite{kupferman-tadmor:euler}).

When applied to equation (\ref{vorticity_CL}), our two-dimensional,
third-order, semi-discrete scheme, (\ref{semi multi_CD})--(\ref{2D_a}), 
takes the form,
\be
\frac{d{\bar \omega}_{j,k}}{dt}=
-\frac{H_{j+1/2,k}^x(t)-H_{j-1/2,k}^x(t)}{\Delta x}
-\frac{H_{j,k+1/2}^y(t)-H_{j,k-1/2}^y(t)}{\Delta y}
+\nu Q_{j,k}(t);
\label{semi inc_Euler}
\ee
with the numerical convection fluxes,
\begin{eqnarray}
H^x_{j+1/2,k}(t) &=&
\frac{u_{j+1/2,k}(t)}{2}\left[{\omega}_{j+1/2,k}^+(t)+
{\omega}_{j+1/2,k}^-(t)\right]- \nonumber \\
&&-\frac{a^x_{j+1/2,k}(t)}{2}\left[{\omega}_{j+1/2,k}^+(t)-
{\omega}_{j+1/2,k}^-(t)\right],\nonumber\\
\label{inc_Euler_fluxes}\\
H^y_{j,k+1/2}(t) &=&
\frac{v_{j,k+1/2}(t)}{2}\left[{\omega}_{j,k+1/2}^+(t)+
{\omega}_{j,k+1/2}^-(t)\right]- \nonumber\\
&&-\frac{a^y_{j,k+1/2}(t)}{2}\left[{\omega}_{j,k+1/2}^+(t)-
{\omega}_{j,k+1/2}^-(t)\right],\nonumber
\end{eqnarray}
and the local speeds,
\be
a^x_{j+1/2,k}(t):=|u_{j+1/2,k}(t)|, \qquad
a^y_{j,k+1/2}(t):=|v_{j,k+1/2}(t)|.
\label{inc_Euler_a}
\ee
To approximate the linear viscosity, $\Delta \omega$, we used
the fourth-order central differencing,
\bea
Q_{j,k}(t)&=&\frac{-{\omega}_{j+2,k}(t)+16{\omega}_{j+1,k}(t)-
30{\omega}_{j,k}(t)+16{\omega}_{j-1,k}(t)-{\omega}_{j-2,k}(t)}{12{\Delta x}^2}+ \nonumber \\
\nonumber\\
&&+\frac{-{\omega}_{j,k+2}(t)+16{\omega}_{j,k+1}(t)-
30{\omega}_{j,k}(t)+16{\omega}_{j,k-1}(t)-{\omega}_{j,k-2}(t)}{12{\Delta y}^2}.
\label{delta omega}
\eea
To compute the intermediate values of the vorticity,
we use the `dimension by dimension' approach described in \S \ref{subsec:2D}:
we reconstruct the corresponding CWENO
interpolants in the $x$- and $y$-directions to obtain the values of
${\omega}_{j+1/2,k}^{\pm}$ and
${\omega}_{j,k+1/2}^{\pm}$.

Another important point in the incompressible computations is that
in every time step one has to recover the velocities,
$\{u_{j,k},v_{j,k}\}$, from the known values of the vorticity,
$\{{\omega}_{j,k}\}$. This can be done in many different ways
(consult, e.g., \cite{levy-tadmor} and the references therein).
Here we have used a stream-function, $\psi$, such that
$\triangle \psi=-\omega$, which is obtained by solving the
nine-points Laplacian, $\triangle {\psi}_{j,k}=-{\omega}_{j,k}(t)$.
This provides the values of the stream-function with 
fourth-order accuracy.   Its gradient, $\nabla \psi$, then recovers
the velocity field,
\bea
u_{j,k}(t)\!\!\!&=&\!\!\!\frac{-{\psi}_{j,k+2}+8{\psi}_{j,k+1}-
8{\psi}_{j,k-1}+{\psi}_{j,k-2}}{12\Delta y},\nonumber\\
\label{velocities}\\
v_{j,k}(t)\!\!\!&=&\!\!\!\frac{{\psi}_{j+2,k}-8{\psi}_{j+1,k}+
8{\psi}_{j-1,k}-{\psi}_{j-2,k}}{12\Delta x}.\nonumber
\eea
\begin{remarks}
\begin{enumerate}
\item Observe that in this way we retain the discrete incompressibility,
namely the discrete velocities computed in (\ref{velocities}) satisfy
\[
\frac{-u_{j+2,k}+8u_{j+1,k}-8u_{j-1,k}+u_{j-2,k}}{12\Delta x}+
\frac{-v_{j,k+2}+8v_{j,k+1}-8v_{j,k-1}+v_{j,k-2}}{12\Delta y}=0.
\]
\item The point-values of the vorticity, which are required for using the nine-points
Laplacian, were computed from its cell averages using the
`dimension by dimension' recipe,
(\ref{dim*dim_1})--(\ref{dim*dim_2}).
\end{enumerate}
\end{remarks}

Finally, the intermediate values of velocities can be computed, e.g.,
using fourth-order averaging,
\bea
u_{j+1/2,k}(t)\!\!\!&=&\!\!\!\frac{-u_{j+2,k}(t)+9u_{j+1,k}(t)+
9u_{j-1,k}(t)-u_{j-2,k}(t)}{16},\nonumber\\
\label{4th-order_aver}\\
v_{j,k+1/2}(t)\!\!\!&=&\!\!\!\frac{-v_{j,k+2}(t)+9v_{j,k+1}(t)+
9v_{j,k-1}(t)-v_{j,k-2}(t)}{16}.\nonumber
\eea

\bigskip
We start our numerical experiments by checking the accuracy of our
scheme, (\ref{semi inc_Euler})--(\ref{4th-order_aver}), augmented with
the DUMKA3 time discretization. We consider the Navier-Stokes equations,
(\ref{inc_Euler})--(\ref{div_free}) with $\nu=0.05$, subject to the
smooth periodic initial data (taken from \cite{Chorin}),
\be
u(x,y,0)=-\cos(x)\sin(y), \qquad v(x,y,0)=\sin(x)\cos(y),
\label{ID Chorin}
\ee
The exact solution to this problem
is simply an exponential decay of the initial data, given by
\[
u(x,y,t)=-\cos(x)\sin(y)e^{-2\nu t}, \qquad
v(x,y,t)=\sin(x)\cos(y)e^{-2\nu t}.
\]

The approximate solution with different number of grid points
was computed at time $t=2$. The errors,
measured in terms of vorticity in
the $L^{\infty}$-, $L^1$- and $L^2$-norms are shown
in Table \ref{table:NS accuracy}. 

\begin{table}[!h]
\begin{center}
\begin{tabular}{|c|cccccc|} \hline
&&&&&&\\
Nx*Ny & $L^{\infty}$-error & rate & $L^1$-error & rate & $L^2$-error & rate \\
&&&&&&\\ \hline
&&&&&&\\
32*32 & 2.429e-02 & -- & 1.791e-01 & -- & 4.559e-02 & -- \\
64*64 & 4.571e-03 & 2.41 & 2.814e-02 & 2.67 & 7.635e-03 & 2.58\\
128*128 & 8.342e-04 & 2.45 & 3.869e-03 & 2.86 & 1.146e-03 & 2.74\\
256*256 & 1.208e-04 & 2.79 & 4.966e-04 & 2.96 & 1.502e-04 & 2.93\\
&&&&&&\\ \hline
\end{tabular}
\caption{Accuracy Test for the Navier-Stokes Equations.
  (\ref{inc_Euler})--(\ref{div_free}), (\ref{ID Chorin}), $\nu=0.05$. Errors at $T=2$
\label{table:NS accuracy}}
\end{center}
\end{table}

Next, the third-order semi-discrete scheme,
(\ref{semi inc_Euler})--(\ref{4th-order_aver}), 
was implemented for the periodic double shear-layer model problem 
taken from \cite{BCG}. First, we solve the Euler equations,
(\ref{inc_Euler})--(\ref{div_free}) with $\nu=0$,
subject to the $(2\pi,2\pi)$-periodic initial data,
\be
u(x,y,0)=\left\{
\begin{array}{ll}
\tanh(\frac{1}{\rho}(y-\pi/2), & y\leq\pi,\\
\\
\tanh(\frac{1}{\rho}(3\pi/2-y), & y>\pi,
\end{array}
\right.\qquad
v(x,y,0)=\delta\cdot\sin(x).
\label{inc_Euler_ID}
\ee
Here, the "thick" shear-layer width parameter, $\rho$, is taken as
$\frac{\pi}{15}$ and the perturbation parameter $\delta=0.05$.

The numerical results at times $T=4, 6, 10$ with $N\!=\!64 \times 64$
and $N\!=\! 128 \times 128$ grid points are presented in
Figures \ref{inc_Euler_contour_64_4}--\ref{inc_Euler_contour_128_10} and
\ref{inc_Euler_64_10}--\ref{inc_Euler_128_10}. 
In order to compare the quality of the results obtained with our new 
method, to previous results, we plot in Figures 
\ref{second-order_contour_64_10}--\ref{second-order_contour_128_10}
the results obtained for the same double shear-layer problem with the
second-order central scheme proposed in \cite{levy-tadmor}.
Compared with the second-order method, the new third-order method, 
can better resolve the large gradients.  Since we are using only
an essentially non-oscillatory reconstruction, some oscillations
are created with the third-order method (and not with the 
``fully'' non-oscillatory second-order method).

Finally, we solve the Navier-Stokes (N-S) equations,
(\ref{inc_Euler})--(\ref{div_free}) with $\nu=0.01$,
augmented with the "thick" shear-layer periodic initial data,
(\ref{inc_Euler_ID}).

The numerical results at time $T=10$ with $N\!=\!64 \times 64$
and $N\!=\! 128 \times 128$ grid points are presented in
Figures \ref{inc_NS_contour_64_10}--\ref{inc_NS_128_10}.  

\begin{center}
\begin{minipage}[t]{7.6cm}
\begin{figure}[H]
\hspace*{-0.5cm}
Contact Author for Figure If Necessary
\caption{Incompressible Euler Equations; Third-order method; {\bf T=4}, 64*64 grid.}
\label{inc_Euler_contour_64_4}
\end{figure}
\end{minipage}\ \hspace{0.4cm}\
\begin{minipage}[t]{7.6cm}
\begin{figure}[H]
Contact Author for Figure If Necessary
\caption{Incompressible Euler Equations; Third-order method; {\bf T=4}, 128*128 grid.}
\label{inc_Euler_contour_128_4}
\end{figure}
\end{minipage}
\end{center}
 
\begin{center}
\begin{minipage}[t]{7.6cm}
\begin{figure}[H]
\hspace*{-0.5cm}
Contact Author for Figure If Necessary
\caption{Incompressible Euler Equations; Third-order method; {\bf T=6}, 64*64 grid.}
\label{inc_Euler_contour_64_6}
\end{figure}
\end{minipage}\ \hspace{0.4cm}\
\begin{minipage}[t]{7.6cm}
\begin{figure}[H]
Contact Author for Figure If Necessary
\caption{Incompressible Euler Equations; Third-order method; {\bf T=6}, 128*128 grid.}
\label{inc_Euler_contour_128_6}
\end{figure}
\end{minipage}
\end{center}

\newpage

\begin{center}
\begin{minipage}[t]{7.6cm}
\begin{figure}[H]
\hspace*{-0.5cm}
Contact Author for Figure If Necessary
\caption{Incompressible Euler Equations; Third-order method; {\bf T=10}, 64*64 grid.}
\label{inc_Euler_contour_64_10}
\end{figure}
\end{minipage}\ \hspace{0.4cm}\
\begin{minipage}[t]{7.6cm}
\begin{figure}[H]
Contact Author for Figure If Necessary
\caption{Incompressible Euler Equations; Third-order method; {\bf T=10}, 128*128 grid.}
\label{inc_Euler_contour_128_10}
\end{figure}
\end{minipage}
\end{center}

\begin{center}
\begin{minipage}[t]{7.6cm}
\begin{figure}[H]
\hspace*{-0.5cm}
Contact Author for Figure If Necessary
\caption{Incompressible Euler Equations; Second-order method; {\bf T=10}, 64*64 grid.}
\label{second-order_contour_64_10}
\end{figure}
\end{minipage}\ \hspace{0.4cm}\
\begin{minipage}[t]{7.6cm}
\begin{figure}[H]
Contact Author for Figure If Necessary
\caption{Incompressible Euler Equations; Second-order method; {\bf T=10}, 128*128 grid.}
\label{second-order_contour_128_10}
\end{figure}
\end{minipage}
\end{center}

\begin{center}
\begin{minipage}[t]{7.6cm}
\begin{figure}[H]
\hspace*{-0.5cm}
Contact Author for Figure If Necessary
\caption{Incompressible Euler Equations; Third-order method; {\bf T=10}, 64*64 grid.}
\label{inc_Euler_64_10}
\end{figure}
\end{minipage}\ \hspace{0.4cm}\
\begin{minipage}[t]{7.6cm}
\begin{figure}[H]
Contact Author for Figure If Necessary
\caption{Incompressible Euler Equations; Third-order method; {\bf T=10}, 128*128 grid.}
\label{inc_Euler_128_10}
\end{figure}
\end{minipage}
\end{center}

\begin{center}
\begin{minipage}[t]{7.6cm}
\begin{figure}[H]
\hspace*{-0.5cm}
Contact Author for Figure If Necessary
\caption{Incompressible Navier-Stokes Equations; Third-order method; {\bf T=10}, 64*64 grid.}
\label{inc_NS_contour_64_10}
\end{figure}
\end{minipage}\ \hspace{0.4cm}\
\begin{minipage}[t]{7.6cm}
\begin{figure}[H]
Contact Author for Figure If Necessary
\caption{Incompressible Navier-Stokes Equations; Third-order method; {\bf T=10}, 128*128 grid.}
\label{inc_NS_contour_128_10}
\end{figure}
\end{minipage}
\end{center}

\begin{center}
\begin{minipage}[t]{7.6cm}
\begin{figure}[H]
\hspace*{-0.5cm}
Contact Author for Figure If Necessary
\caption{Incompressible Navier-Stokes Equations; Third-order method; {\bf T=10}, 64*64 grid.}
\label{inc_NS_64_10}
\end{figure}
\end{minipage}\ \hspace{0.4cm}\
\begin{minipage}[t]{7.6cm}
\begin{figure}[H]
Contact Author for Figure If Necessary
\caption{Incompressible Navier-Stokes Equations; Third-order method; {\bf T=10}, 128*128 grid.}
\label{inc_NS_128_10}
\end{figure}
\end{minipage}
\end{center}

\begin{acknowledgment}
The authors would like to thank Prof.\ S.\ Karni and Prof.\ R.\ Krasny for helpful comments.
The work of A.\ Kurganov was supported in part by the NSF Group Infrastructure Grant.
The work of D.\ Levy was supported in part by the Applied Mathematical Sciences 
subprogram of the Office of Science, U.S.\ Department of Energy, under 
contract DE--AC03--76--SF00098.
Part of this work was done while A.K.\ was visiting the Lawrence Berkeley Lab.
\end{acknowledgment}



\begin{thebibliography}{99}
\bibitem{arminjon:nt.french} Arminjon P., Viallon M.-C., 
{\em G\'en\'eralisation du Sch\'ema de Nessyahu-Tadmor pour Une 
\'Equation Hyperbolique \`a Deux Dimensions D'espace}, 
C.R. Acad. Sci. Paris, {\bf t. 320} , s\'erie I. (1995), pp.85--88.

\bibitem{BCG} Bell J. B., Colella P., Glaz H. M.,
{\em A Second-Order Projection Method for the Incompressible
Navier-Stokes Equations},
JCP, {\bf 85}, (1989), pp.257--283.

\bibitem{bianco-puppo-russo:central} Bianco F., Puppo G., Russo G., 
{\em High Order Central Schemes for Hyperbolic Systems of Conservation Laws}, 
SIAM J. Sci. Comp., to appear.

\bibitem{Chorin} Chorin A.,
{\em Numerical Solution of the Navier-Stokes Equations},
Math. Comp., {\bf 22}, (1968), pp.745--762.

\bibitem{lxf} Friedrichs K. O., Lax P. D., 
{\em Systems of Conservation Equations with a Convex Extension}, 
Proc. Nat. Acad. Sci., {\bf 68}, (1971), pp.1686--1688.

\bibitem{god-rav:difference} Godlewski E., Raviart P.-A., 
{\em Numerical Approximation of Hyperbolic Systems of Conservation Laws}, 
Springer, New York, 1996.

\bibitem{God}
Godunov S. K., {\em A Finite Difference Method for the Numerical
Computation of Discontinuous Solutions of the Equations of Fluid Dynamics},
Mat. Sb., {\bf 47}, (1959), pp.271--290.

\bibitem{harten:ACM} Harten A., {\em The Artificial Compression
Method for Computation of Shocks and Contact Discontinuities, III.
Self-Adjusting Hybrid Schemes},
Math. Comp., {\bf 32}, (1978), pp.363--389.

\bibitem{harten-eoc:eno} Harten A., Engquist B., Osher S., Chakravarthy S.,
{\em Uniformly High Order Accurate Essentially Non-oscillatory Schemes III}, 
JCP, {\bf 71}, (1987), pp.231--303.

\bibitem{jiang-levy-osher-tadmor:stg} Jiang G.-S., Levy D., Lin C.-T.,
Osher S., Tadmor E. 
{\em High-Resolution Non-Oscillatory Central Schemes with Non-Staggered 
Grids for Hyperbolic Conservation Laws},
SINUM, {\bf 35}, (1998), pp.2147--2168.

\bibitem{jiang-shu:weno} Jiang G.-S., Shu C.-W., 
{\em Efficient Implementation of Weighted ENO Schemes}, 
JCP, {\bf 126}, (1996), pp.202--228.

\bibitem{jiang-tadmor:nonosc} Jiang G.-S., Tadmor E., 
{\em Nonoscillatory Central Schemes for
Multidimensional Hyperbolic Conservation Laws}, 
SIAM J. Sci. Comp., {\bf 19}, (1998), pp.1892--1917. 

\bibitem{kupferman-tadmor:euler} Kupferman R., Tadmor E.,
 {\em A Fast High-Resolution Second-Order Central Scheme for Incompressible Flows}, 
Proc. Nat. Acad. Sci., {\bf 94}, (1997), pp. 4848-4852.

\bibitem{KLR} Kurganov A., Levy D., Rosenau P.,
{\em On Burgers-Type Equations with Nonmonotonic Dissipative Fluxes},
Comm. Pure Appl. Math., {\bf 51}, (1998), pp.443--473.

\bibitem{KR} Kurganov A., Rosenau P.,
{\em Effects of a Saturating Dissipation in Burgers-Type Equations},
Comm. Pure Appl. Math., {\bf 50}, (1997), pp.753--771.

\bibitem{kurganov-tadmor:semi} Kurganov A., Tadmor E.,
{\em New High-Resolution Central Schemes for Nonlinear
 Conservation Laws and Convection-Diffusion Equations}, 
submitted.

\bibitem{vLeV} van Leer B.,
{\em Towards the ultimate conservative difference scheme, V.
A second order sequel to Godunov's method},
JCP, {\bf 32}, (1979), pp. 101--136.

\bibitem{levy:third}  Levy D., 
{\em A Third-order 2D Central Schemes for Conservation Laws}, 
INRIA School on Hyperbolic Systems, Vol. I (1998), pp.489--504.

\bibitem{levy-puppo-russo:1d} Levy D., Puppo G., Russo G., 
{\em Central WENO Schemes for Hyperbolic Systems of Conservation Laws},
M2AN, in press.

\bibitem{levy-puppo-russo:2d-3} Levy D., Puppo G., Russo G., 
{\em A Third Order Central WENO Scheme for 2D Conservation Laws},
Appl. Nume. Math., in press.

\bibitem{levy-puppo-russo:balance} Levy D., Puppo G., Russo G., 
{\em Compact Central WENO Schemes for Multidimensional Conservation Laws},
submitted.

\bibitem{levy-puppo-russo:2d-4} Levy D., Puppo G., Russo G.,
{\em Central WENO Schemes for Multi-Dimensional Hyperbolic Systems of 
 Conservation Laws}, 
in preparation.

\bibitem{levy-puppo-russo:tv} Levy D., Puppo G., Russo G., 
{\em On the Behavior of the Total Variation in CWENO Methods for Conservation Laws}, 
Appl. Nume. Math., in press.

\bibitem{levy-tadmor} Levy D., Tadmor E.,
{\em Non-Oscillatory Central Schemes for the Incompressible 2-D Euler Equations},
Math. Res. Lett., {\bf 4}, (1997), pp.1--20.

\bibitem{liu-osher:nonosc} Liu X.-D., Osher S., 
{\em Nonoscillatory High Order Accurate Self-Similar Maximum Principle 
Satisfying Shock Capturing Schemes I}, 
SINUM, {\bf 33}, no. 2 (1996), pp.760--779.

\bibitem{liu-osher-chan:weno} Liu X.-D., Osher S., Chan T., 
{\em Weighted Essentially Non-oscillatory Schemes}, 
JCP, {\bf 115}, (1994), pp.200--212.

\bibitem{liu-tadmor:3rd} Liu X.-D., Tadmor E., 
{\em Third Order Nonoscillatory Central Scheme
for Hyperbolic Conservation Laws}, 
Numer. Math., {\bf 79}, (1998), pp.397--425. 

\bibitem{Med} Medovikov A.A., {\em High Order Explicit Methods for
Parabolic Equations},
BIT, {\bf 38}, (1998) 2, pp.372--390.

\bibitem{nessyahu-tadmor:non-oscillatory} Nessyahu H., Tadmor E., 
{\em Non-oscillatory Central Differencing for Hyperbolic Conservation Laws}, 
JCP, {\bf 87}, no. 2 (1990), pp.408--463.

\bibitem{roe:approximate} Roe P. L., 
{\em Approximate Riemann Solvers, Parameter Vectors, and Difference Schemes}, 
JCP, {\bf 43}, (1981), pp.357--372.

\bibitem{shu:eno} Shu C.-W., 
{\em Numerical Experiments on the Accuracy of ENO and Modified ENO Schemes}, 
J. Sci. Comp., {\bf 5}, vol. 2, (1990), pp.127--149.

\bibitem{SO} Shu C.-W., Osher S.,
{\em Efficient Implementation of Essentially Non-Oscillatory
Shock-Capturing Schemes},
JCP, {\bf 77}, (1988), pp.439--471.

\bibitem{sod:hyperbolic} Sod G., 
{\em A Survey of Several Finite Difference Methods for Systems of Nonlinear 
Hyperbolic Conservation Laws}, 
JCP, {\bf 22}, (1978), pp.1--31.

\bibitem{tadmor:approximate} Tadmor E., 
{\em Approximate Solutions of Nonlinear Conservation Laws},
CIME Lecture notes, 1997, UCLA CAM Report 97-51.

\end{thebibliography}
\end{document}